\documentclass[reqno]{amsart}

\usepackage{epsfig,graphicx,amsfonts,amssymb,amsmath,amscd,amsthm,color}
\usepackage[all]{xy}
\usepackage{pgf,tikz}
\usetikzlibrary{arrows}

\numberwithin{equation}{section}

\newtheorem{theorem}{Theorem}[section]

\newcommand{\Ker}{{\rm Ker}}
\newcommand{\rank}{{\rm rank}}

\newcommand{\CCC}{{\mathbb C}}

\newcommand{\RR}{{\mathbb R}}
\newcommand{\ZZ}{{\mathbb Z}}
\newcommand{\PP}{{\mathbb P}}

\newcommand{\HH}{{\mathbb H}}
\newcommand{\AAA}{{\mathbb A}}

\newcommand{\CL}{{\mathcal L}}

\newcommand{\CO}{{\mathcal O}}

\newcommand{\CJ}{{\mathcal J}}

\newcommand{\CE}{{\mathcal E}}

\def\:{\colon}

\begin{document}

\title{An example of a reducible Severi variety}
\author{Ilya Tyomkin}

\thanks{The research leading to these results has received funding from the European Union Seventh Framework Programme (FP7/2007-2013) under grant agreement 248826.}

\address{Department of Mathematics, Ben-Gurion University of the Negev, P.O.Box 653, Be'er Sheva, 84105, Israel}
\email{tyomkin@cs.bgu.ac.il}

\begin{abstract}
We construct a positive-dimensional, reducible Severi variety on a toric surface.
\end{abstract}
\maketitle

\section{Introduction}

In this note we discuss the irreducibility problem for Severi varieties, which parameterize irreducible reduced curves of geometric genus $g$ in a given linear system $|\CL|$ on a given projective surface $X$. Throughout the paper we work over the field of complex numbers.

Severi varieties were introduced by Severi \cite{Sev21} in 1920s in an attempt to prove the irreducibility of moduli spaces of algebraic curves of a given genus. He considered the case of plane curves and proved the irreducibility, but the proof contained a gap, and the first complete proof of the irreducibility was obtained by Harris \cite{Har86} only in 1986. The result of Harris was generalized to other rational surfaces. In particular, the irreducibility is known for Hirzebruch surfaces \cite{T07}, and in the case of rational curves it is also known for all toric surfaces \cite{T07} and all del Pezzo surfaces \cite{Tes09} with one trivial zero-dimensional exception: the variety of rational plane cubics through eight points in general position. The goal of this note is to describe an example of a toric surface admitting reducible positive-dimensional Severi varieties.

The idea of Harris's proof of irreducibility of Severi varieties is as follows: A standard deformation-theoretic argument shows that in a neighborhood of an irreducible rational nodal curve the Severi variety consists of several smooth branches parameterized by subsets of $p_a(\CL)-g$ nodes (one branch for each choice of $p_a(\CL)-g$ nodes). So, first, one proves that the monodromy acts as a full symmetric group on the set of nodes, and concludes that there exists a unique component of the Severi variety whose closure in the linear system contains an irreducible rational nodal curve. Second, one proves that the closure of any component of the Severi variety contains a nodal curve of smaller genus, and deduces from this that it contains a rational nodal curve, which completes the proof.

The first step is easy in the plane case and in the case of Hirzebruch surfaces, and the main difficulty of the proof is in the second step, which requires several brilliant ideas and involved arguments from the deformation theory and the theory of curves. As we will see below, in our example the first step is the one that fails, while the second step generalizes smoothly.

\section{The example}

We refer the reader to \cite{Dan78, Ful93} for an introduction to toric geometry. Let $N=\ZZ^2=\ZZ e_1\oplus\ZZ e_2$ be a lattice of rank two, and $M$ the corresponding dual lattice. Consider the complete fan $\Sigma$ in $\RR^2$ with rays generated by $\pm 2e_1\pm e_2$, and let $X$ be the toric surface corresponding to $\Sigma$.
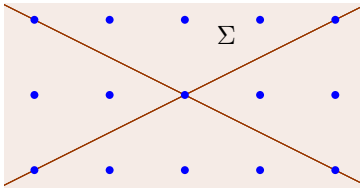
\begin{center}
\definecolor{qqqqff}{rgb}{0,0,1}
\definecolor{zzttqq}{rgb}{0.6,0.2,0}
\definecolor{zzwwqq}{rgb}{0.6,0.4,0}
\begin{tikzpicture}[line cap=round,line join=round,>=triangle 45,x=1.0cm,y=1.0cm]
\clip(-2.4,1.74) rectangle (2.32,4.22);
\fill[color=zzttqq,fill=zzttqq,fill opacity=0.1] (0,3) -- (4,5) -- (-4,5) -- (4,1) -- (-4,1) -- cycle;
\fill[color=zzttqq,fill=zzttqq,fill opacity=0.1] (0,3) -- (-4,5) -- (-4,1) -- (4,5) -- (4,1) -- cycle;
\draw [color=zzwwqq,domain=0.0:2.320000000000002] plot(\x,{(--6--1*\x)/2});
\draw [color=zzwwqq,domain=0.0:2.320000000000002] plot(\x,{(--6-1*\x)/2});
\draw [color=zzwwqq,domain=-2.400000000000003:0.0] plot(\x,{(-6-1*\x)/-2});
\draw [color=zzwwqq,domain=-2.400000000000003:0.0] plot(\x,{(-6--1*\x)/-2});
\draw [color=zzttqq] (0,3)-- (4,5);
\draw [color=zzttqq] (4,5)-- (-4,5);
\draw [color=zzttqq] (-4,5)-- (4,1);
\draw [color=zzttqq] (4,1)-- (-4,1);
\draw [color=zzttqq] (-4,1)-- (0,3);
\draw [color=zzttqq] (0,3)-- (-4,5);
\draw [color=zzttqq] (-4,5)-- (-4,1);
\draw [color=zzttqq] (-4,1)-- (4,5);
\draw [color=zzttqq] (4,5)-- (4,1);
\draw [color=zzttqq] (4,1)-- (0,3);
\draw (0.3,4.04) node[anchor=north west] {$\Sigma$};
\begin{scriptsize}
\fill [color=qqqqff] (-2,4) circle (1.5pt);
\fill [color=qqqqff] (-1,4) circle (1.5pt);
\fill [color=qqqqff] (0,4) circle (1.5pt);
\fill [color=qqqqff] (1,4) circle (1.5pt);
\fill [color=qqqqff] (2,4) circle (1.5pt);
\fill [color=qqqqff] (2,3) circle (1.5pt);
\fill [color=qqqqff] (2,2) circle (1.5pt);
\fill [color=qqqqff] (1,2) circle (1.5pt);
\fill [color=qqqqff] (1,3) circle (1.5pt);
\fill [color=qqqqff] (0,3) circle (1.5pt);
\fill [color=qqqqff] (0,2) circle (1.5pt);
\fill [color=qqqqff] (-1,2) circle (1.5pt);
\fill [color=qqqqff] (-1,3) circle (1.5pt);
\fill [color=qqqqff] (-2,3) circle (1.5pt);
\fill [color=qqqqff] (-2,2) circle (1.5pt);
\end{scriptsize}
\end{tikzpicture}
\end{center}

Let $\Delta$ be the lattice polygon in $M_\RR:=M\otimes_\ZZ\RR$ with vertices $\pm e^1, \pm 2e^2$, where $\{e^i\}\subset M$ denotes the dual basis to the basis $\{e_i\}\subset N$. It defines an ample line bundle on the surface $X$, denoted by $\CL$. Since the arithmetic genus of a line bundle is equal to the number of inner lattice points in the corresponding polygon we have: $p_a(\CL)=3$.
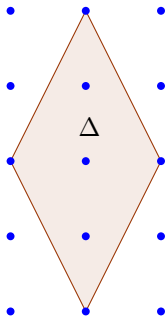
\begin{center}
\definecolor{zzttqq}{rgb}{0.6,0.2,0}
\definecolor{qqqqff}{rgb}{0,0,1}
\begin{tikzpicture}[line cap=round,line join=round,>=triangle 45,x=1.0cm,y=1.0cm]
\clip(0.76,-0.3) rectangle (3.22,4.24);
\fill[color=zzttqq,fill=zzttqq,fill opacity=0.1] (2,4) -- (1,2) -- (2,0) -- (3,2) -- cycle;
\draw [color=zzttqq] (2,4)-- (1,2);
\draw [color=zzttqq] (1,2)-- (2,0);
\draw [color=zzttqq] (2,0)-- (3,2);
\draw [color=zzttqq] (3,2)-- (2,4);
\draw (9.72,2.04) node[anchor=north west] {$\Delta_0$};
\draw (1.78,2.7) node[anchor=north west] {$\Delta$};
\begin{scriptsize}
\fill [color=qqqqff] (2,2) circle (1.5pt);
\fill [color=qqqqff] (1,2) circle (1.5pt);
\fill [color=qqqqff] (1,3) circle (1.5pt);
\fill [color=qqqqff] (2,3) circle (1.5pt);
\fill [color=qqqqff] (3,3) circle (1.5pt);
\fill [color=qqqqff] (3,2) circle (1.5pt);
\fill [color=qqqqff] (1,1) circle (1.5pt);
\fill [color=qqqqff] (2,1) circle (1.5pt);
\fill [color=qqqqff] (3,1) circle (1.5pt);
\fill [color=qqqqff] (2,0) circle (1.5pt);
\fill [color=qqqqff] (3,0) circle (1.5pt);
\fill [color=qqqqff] (1,0) circle (1.5pt);
\fill [color=qqqqff] (1,4) circle (1.5pt);
\fill [color=qqqqff] (2,4) circle (1.5pt);
\fill [color=qqqqff] (3,4) circle (1.5pt);
\end{scriptsize}
\end{tikzpicture}
\end{center}
For the convenience of a reader not familiar with toric geometry, we provide an explicit description of the immersion of the surface $X$ into the projective space $\PP^6$ given by the linear system $|\CL|$ at the very end of \S~\ref{sec:Ran}.

The Severi variety we are going to consider is the variety parameterizing irreducible reduced curves of genus one in the linear system $|\CL|$ on the surface $X$ that do not pass through the zero-dimensional orbits. We denote this variety by $V:=V^{irr}(X,\CL,1)$.

\begin{theorem}\label{thm:main} Severi variety $V$ is reducible. It is equidimensional of dimension $4$, has exactly two irreducible components, and general points of each component correspond to irreducible reduced nodal curves of genus one. Furthermore, the closure of each component in $|\CL|$ contains points corresponding to irreducible rational nodal curves.
\end{theorem}

\section{Proof of Theorem~\ref{thm:main}}
\subsection{Curves in $|\CL|$}
Let $C\in |\CL|$ be an irreducible curve that does not contain the zero-dimensional orbits of $X$, and $E\to C$ the normalization of $C$. Denote by $f$ the composition morphism $E\to C\hookrightarrow X$. Let $n_1:=2e_1+e_2$, $n_2:=-2e_1+e_2$, $n_3:=-2e_1-e_2$, and $n_4:=2e_1-e_2$ be the $N$-generators of the rays of $\Sigma$, and $D_1,\dotsc,D_4\subset X$ the corresponding divisors. Since $C.D_i=1$ for all $i$, there exist points $p_1,\dotsc,p_4\in E$ such that $f^*D_i=p_i$. Furthermore, we obtain a homomorphism $f^*\:M\to \CO_E(E\setminus\cup_ip_i)\subset \CCC(E)$ given by $m\mapsto f^*(x^m)$, which satisfies the following:
\begin{equation}\label{eq:div}
div(f^*(x^m))=\sum_{i=1}^4(m,n_i)p_i
\end{equation}
for any $m\in M$. Vice versa, given four distinct points $p_1,\dotsc,p_4\in E$ and a homomorphism $\phi\:M\to \CO_E(E\setminus\cup_ip_i)\subset \CCC(E)$ satisfying \eqref{eq:div} for any $m\in M$, we obtain a morphism $f\:E\setminus\cup_ip_i\to T_N\subset X$ that extends to $E$ such that $f(p_i)\in D_i$ for all $i$, and $\phi=f^*$. Note that to give a homomorphism $\phi\:M\to \CO_E(E\setminus\cup_ip_i)\subset \CCC(E)$ satisfying \eqref{eq:div} for any $m\in M$ it is sufficient to specify the values of $\phi$ on any set of generators of $M$ and to verify \eqref{eq:div} on these generators. Plainly, the data $(E;p_1,\dotsc,p_4;\phi)$ and $(E';p'_1,\dotsc,p'_4;\phi')$ correspond to the same curve $C$ if and only if there exists an isomorphism $\alpha\:E\to E'$ such that $\alpha(p_i)=p'_i$ for all $i$, and $\phi=\alpha^*\circ\phi'$.

\subsection{The two components of $V$}
Let $(E;p_1,\dotsc,p_4;\phi)$ be a datum as above such that the genus of $E$ is one. Then $(E,p_1)$ is an elliptic curve, and \eqref{eq:div} for the generators $\{e^1+2e^2, e^2\}$ of $M$ is equivalent to the following identities in the group $E$: $4[p_3]=0$ and $[p_2]=[p_3]+[p_4]$. We conclude that the Severi variety $V$ is dominated by the variety $W$ parameterizing the following data: $(E,O;p,q;x_1,x_2)$, where $(E,O)$ is an elliptic curve, $p\in E$ is a point of order $2$ or $4$, $q\in E$ is any point such that $[q]\ne 0, \pm [p]$, and $x_1, x_2$ are rational functions with divisors $4O-4p$ and $p+q-r-O$ respectively, where $[r]=[p]+[q]$. It also follows from the description above that the projection $W\to V$ is at most two-to-one.

We see that $W$ consists of two components $W_2$ and $W_4$ distinguished by the order of the point $p$, and hence the projections of $W_2$ and $W_4$ to $V$ are disjoint. By forgetting the rational functions $x_1, x_2$ one represents $W_k$ as a $\CCC^*\times\CCC^*$-torsor over a dense open subset of the universal curve $\CE_k$ over the moduli space $X_1[k]$ of elliptic curve with the (partial) level-$k$ structure. Recall that $X_1[k]\simeq\Gamma_1[k]\backslash\HH$, where $\HH$ is the upper half space and $\Gamma_1[k]$ the congruence subgroup
$$\Gamma_1[k]=\left\{\left.
A\in SL(2,\ZZ)\,\right| A\equiv
\left(
        \begin{array}{cc}
             1 & * \\
             0 & 1 \\
        \end{array}
    \right) \mod k
\right\}.$$
Thus, $W_k$ are irreducible since so are the moduli spaces $X_1[k]$ and the universal curves $\CE_k$. Furthermore, $\dim(W_k)=\dim(\CCC^*\times\CCC^*)+\dim(\CE_k)=4$. Thus, $V$ has exactly two irreducible components, each of dimension $4$. The claim about the nodality of the general curve parameterized by each irreducible component follows from \cite[Theorem~2.8]{T07} applied to a toric desingularization of $X$.

\subsection{Irreducible rational nodal curves corresponding to points in the closure of the components}
For $k=1,2$, consider the maps $f_k\:\PP^1\to X$ given by $f_k^*(x^m)=t^{(m,n_k)}$, where $t$ is the coordinate on $\AAA^1\subset\PP^1$. Set $C_k:=f_k(\PP^1)$, $C:=C_1\cup C_2$, and $p_i:=f^{-1}(D_i)$ for $i=1,\dotsc,4$. Then $C$ is a $4$-nodal curve in the linear system $|\CL|$. Let $q\in C$ be one of its nodes. We claim that $C$, and any small deformation $C_\epsilon$ of $C$ that smooths out $q$, belong to the closures $\overline{W}_i$, for $i=2,4$.

First, a standard deformation-theoretic argument shows that the nodes of $C$ can be smoothed out independently. Indeed, the tangent space to the space of deformation preserving a given collection of nodes $\{r_j\}$ is given by $H^0(C, \CJ\otimes\CL|_C)$, where $\CJ$ denotes the ideal of the subscheme $\{r_j\}\subset C$. It follows from the Riemann-Roch theorem that $h^0(C, \CJ\otimes\CL|_C)=h^0(C,\CL|_C)-|\{r_j\}|$, and hence the nodes impose independent conditions on the deformations (see the statements and the proofs of \cite[Proposition~2.11 and Claim~2.12]{T07} for details). Thus, $C$ and $C_\epsilon$ belong to the closure of $V$. To prove that they belong to the closure of $W_i$ it remains to show that we can choose a node $r\in C$ different from $q$ such that the order of $f^*(\CO_X(D_3-D_1))$ in the Jacobian $J(E)$ is $i$, where $f\:E\to X$ denotes the partial normalization of $C$ preserving the nodes $q,r$.

Recall that there exists a natural exact sequence (cf. \cite[p.89]{ACG11})
$$1\to \CCC^\times\to J(E)\to J(C_1)\times J(C_2)\to 0,$$
and since $J(C_i)=J(\PP^1)=0$, we have $J(E)\simeq\CCC^\times$. Furthermore, the class of a bi-degree $(0,0)$ divisor $\CO_E(D)\in J(E)\simeq\CCC^\times$ disjoint from the nodes of $E$ can be computed as $s_1(a)s_2(b)/s_1(b)s_2(a)$, where $s_k\in \CCC(C_k)$ are rational functions for which $div(s_k)=D\cap C_k$, and $a,b\in E$ are the nodes. Thus, the class of the divisor $f^*(\CO_X(D_3-D_1))$ is the ratio $t_1(a)/t_1(b)$, where $t_1$ is the coordinate on $C_1$ vanishing at $p_1$ and having a pole at $p_3$. So, let us compute the $t_1$-coordinates of the nodes of $C$. By the definition of $C$ and of $f_i$, they satisfy the equation $t_1^{(m,n_1)}=t_2^{(m,n_2)}$ for all $M$. The latter is equivalent to $t_1=t_2$ and $t_1^4=1$, i.e., $t_1$-coordinates of the nodes are precisely the roots of unity of order four. Plainly, the order of the ratio $t_1(a)/t_1(b)$ is equal to the order of $t_1(r)/t_1(q)$, and hence the order of $f^*(\CO_X(D_3-D_1))\in J(E)$ can be two or four depending on the choice of $r$.

\section{Final remarks}

\subsection{Monodromy action on the nodes of rational curves}
Any irreducible rational curve $D$ is a deformation $C_\epsilon$ of the curve $C$ from the last part of the proof of the Theorem. It is not difficult to see that the monodromy group does not act transitively on the nodes of $C_\epsilon$. In fact it acts as a transposition on a pair of nodes and preserves the third node. This is precisely the point where the generalization of Harris's proof fails.

\subsection{Toric description of the components}
One can distinguish between the irreducible components geometrically as follows: Let $N'\subset N$ be the sublattice generated by $2e_1$ and $e_2$, and $M\subset M'$ the corresponding dual lattice. Let $X'$ be the toric surfaces corresponding to $\Sigma$ with respect to the lattice $N'$. Then there exists a natural projection $\pi\:X'\to X$ that corresponds to the embedding of the lattices $N'\hookrightarrow N$, and $X$ is the quotient of $X'$ by the action of the group $G:=\Ker(T_N\to T_{N'})\simeq \mu_2$. Furthermore, the action of $G$ on $X'$ is free away from the zero-dimensional orbits.
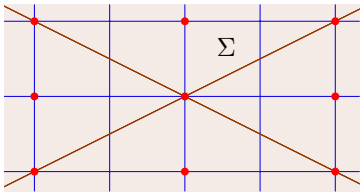
\begin{center}
\definecolor{zzttqq}{rgb}{0.6,0.2,0}
\definecolor{zzwwqq}{rgb}{0.6,0.4,0}
\definecolor{ffqqqq}{rgb}{1,0,0}
\definecolor{qqqqff}{rgb}{0,0,1}
\begin{tikzpicture}[line cap=round,line join=round,>=triangle 45,x=1.0cm,y=1.0cm]
\draw [color=qqqqff,, xstep=1.0cm,ystep=1.0cm] (-2.4,1.74) grid (2.32,4.22);
\clip(-2.4,1.74) rectangle (2.32,4.22);
\fill[color=zzttqq,fill=zzttqq,fill opacity=0.1] (0,3) -- (4,5) -- (-4,5) -- (4,1) -- (-4,1) -- cycle;
\fill[color=zzttqq,fill=zzttqq,fill opacity=0.1] (0,3) -- (-4,5) -- (-4,1) -- (4,5) -- (4,1) -- cycle;
\draw [color=zzwwqq,domain=0.0:2.320000000000002] plot(\x,{(--6--1*\x)/2});
\draw [color=zzwwqq,domain=0.0:2.320000000000002] plot(\x,{(--6-1*\x)/2});
\draw [color=zzwwqq,domain=-2.400000000000003:0.0] plot(\x,{(-6-1*\x)/-2});
\draw [color=zzwwqq,domain=-2.400000000000003:0.0] plot(\x,{(-6--1*\x)/-2});
\draw [color=zzttqq] (0,3)-- (4,5);
\draw [color=zzttqq] (4,5)-- (-4,5);
\draw [color=zzttqq] (-4,5)-- (4,1);
\draw [color=zzttqq] (4,1)-- (-4,1);
\draw [color=zzttqq] (-4,1)-- (0,3);
\draw [color=zzttqq] (0,3)-- (-4,5);
\draw [color=zzttqq] (-4,5)-- (-4,1);
\draw [color=zzttqq] (-4,1)-- (4,5);
\draw [color=zzttqq] (4,5)-- (4,1);
\draw [color=zzttqq] (4,1)-- (0,3);
\draw (0.3,3.9) node[anchor=north west] {$\Sigma$};
\begin{scriptsize}
\fill [color=ffqqqq] (-2,4) circle (1.5pt);
\fill [color=ffqqqq] (0,4) circle (1.5pt);
\fill [color=ffqqqq] (2,4) circle (1.5pt);
\fill [color=ffqqqq] (-2,3) circle (1.5pt);
\fill [color=ffqqqq] (0,3) circle (1.5pt);
\fill [color=ffqqqq] (2,3) circle (1.5pt);
\fill [color=ffqqqq] (-2,2) circle (1.5pt);
\fill [color=ffqqqq] (0,2) circle (1.5pt);
\fill [color=ffqqqq] (2,2) circle (1.5pt);
\end{scriptsize}
\end{tikzpicture}
\end{center}
In the picture above the lattice $N$ is blue, and $N'$ is red. The polygon $\Delta$ is integral with respect to $M'$, and hence defines an ample line bundle $\CL'$ on $X'$ which satisfies: $\pi^*\CL=\CL'$. Let $\Delta_0\subset M'_\RR$ be the lattice polygon with vertices $\pm \frac{1}{2}e^1, \pm e^2$, and $\CL_0$ the corresponding line bundle on $X'$. Then $\CL'=\CL_0^{\otimes 2}$. Furthermore, $p_a(\CL')=5$ and $p_a(\CL_0)=1$.
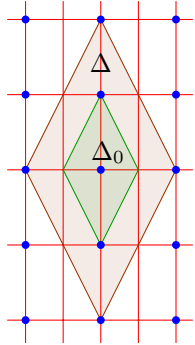
\begin{center}
\definecolor{qqzzqq}{rgb}{0,0.6,0}
\definecolor{zzttqq}{rgb}{0.6,0.2,0}
\definecolor{qqqqff}{rgb}{0,0,1}
\definecolor{ffqqtt}{rgb}{1,0,0}
\begin{tikzpicture}[line cap=round,line join=round,>=triangle 45,x=1.0cm,y=1.0cm]
\draw [color=ffqqtt,, xstep=0.5cm,ystep=1.0cm] (0.76,-0.3) grid (3.22,4.24);
\clip(0.76,-0.3) rectangle (3.22,4.24);
\fill[color=zzttqq,fill=zzttqq,fill opacity=0.1] (2,4) -- (1,2) -- (2,0) -- (3,2) -- cycle;
\fill[color=qqzzqq,fill=qqzzqq,fill opacity=0.1] (2,3) -- (1.5,2) -- (2,1) -- (2.5,2) -- cycle;
\draw [color=zzttqq] (2,4)-- (1,2);
\draw [color=zzttqq] (1,2)-- (2,0);
\draw [color=zzttqq] (2,0)-- (3,2);
\draw [color=zzttqq] (3,2)-- (2,4);
\draw [color=qqzzqq] (2,3)-- (1.5,2);
\draw [color=qqzzqq] (1.5,2)-- (2,1);
\draw [color=qqzzqq] (2,1)-- (2.5,2);
\draw [color=qqzzqq] (2.5,2)-- (2,3);
\draw (1.75,2.5) node[anchor=north west] {$\Delta_0$};
\draw (1.73,3.66) node[anchor=north west] {$\Delta$};
\begin{scriptsize}
\fill [color=qqqqff] (2,2) circle (1.5pt);
\fill [color=qqqqff] (1,2) circle (1.5pt);
\fill [color=qqqqff] (1,3) circle (1.5pt);
\fill [color=qqqqff] (2,3) circle (1.5pt);
\fill [color=qqqqff] (3,3) circle (1.5pt);
\fill [color=qqqqff] (3,2) circle (1.5pt);
\fill [color=qqqqff] (1,1) circle (1.5pt);
\fill [color=qqqqff] (2,1) circle (1.5pt);
\fill [color=qqqqff] (3,1) circle (1.5pt);
\fill [color=qqqqff] (2,0) circle (1.5pt);
\fill [color=qqqqff] (3,0) circle (1.5pt);
\fill [color=qqqqff] (1,0) circle (1.5pt);
\fill [color=qqqqff] (1,4) circle (1.5pt);
\fill [color=qqqqff] (2,4) circle (1.5pt);
\fill [color=qqqqff] (3,4) circle (1.5pt);
\end{scriptsize}
\end{tikzpicture}
\end{center}
In the picture above the lattice $M$ is blue, and $M'$ is red.

Let $C\in V$ be any curve, and $E\to C$ its normalization. Then $E\times_XX'\to E$ is an unramified covering of degree two, and hence $E\times_XX'$ is either a curve of genus one or a disjoint union of two such curves by the Riemann-Hurwitz formula. Hence the preimage $\pi^{-1}(C)\in |\CL'|$ is either a $\mu_2$-equivariant element of the Severi variety $V^{irr}(X',\CL',1)$ or a union of two irreducible curves $C_1\cup C_2$, where $C_1\in V^{irr}(X',\CL_0,1)=|\CL_0|$ is arbitrary, $C_2=\sigma(C_1)$, and $\sigma\in\mu_2$ is the non-identity element. The locus of curves $C$ of first type is precisely $W_4$, and of second type is $W_2$. Notice, that this description provides an explicit double covering of $W_2$ by an open subset of the projective $4$-space $|\CL_0|$.

\subsection{The degrees of the components $W_2$ and $W_4$}
One can use Mikhalkin's correspondence theorem \cite{Mikh05} (see also \cite{Tyo10}), to compute the degree of $V$. It turns out that $\deg(V)=34$. The description of $W_2$ from the previous paragraph allows one to compute the degree of $W_2$ easily. Indeed, pick four points $q_1,\dotsc, q_4\in T_N$ in general position. Then the degree of $W_2$ is equal to the number of curves in $W_2$ passing through $q_1,\dotsc, q_4$. To find this number, we observe that for each lifting of $q_1,\dotsc, q_4$ there exists a unique curve in $|\CL_0|$ passing through the liftings. Thus, the number of curves in $\pi_*|\CL_0|$ through $q_1,\dotsc, q_4$ is $\frac{1}{2}2^4=8$ since each such curve is the image of exactly two curves in $|\CL_0|$. Thus, $\deg(W_2)=8$, and hence $\deg(W_4)=26$.

It is not difficult to check that the non-archimedean amoeba of a curve $C\in |\CL_0|$ coincides with the amoeba of $\pi(C)$. In particular, the amoebas of curves parameterized by $W_2$ are dual to $M'$-subdivisions of the polygon $\Delta_0$. However, by Mikhalkin's realization theorem, there exist curves of genus one in $|\CL|$ whose amoebas are not dual to any $M'$-subdivision of $\Delta_0$, e.g., amoeba dual to the subdivision of $\Delta$ into four triangles with vertices $\pm e^1, \pm2e^2, e^2$. One can use this observation or the computation of the degrees of $V$ and $W_2$ to prove {\em tropically} that $V$ has at least two irreducible components!

\subsection{Computational aspects}
The Severi variety $V$ has codimension two in the linear system $|\CL|$. It is the complement of the loci $U_c, U_1,\dotsc, U_4$ in the singular locus of the discriminant of $|\CL|$, where $U_c$ denotes the locus of cuspidal curves and $U_1,\dotsc, U_4$ are the loci of curves passing through the singular points of $X$. Bernd Sturmfels and Diane Maclagan used Macaulay2 to compute explicitly the ideals of the components $W_2$ and $W_4$, their degrees, etc.

\subsection{Generalizations}\label{sec:Ran}

The example presented in this paper can be generalized to a certain class of toric surfaces, and to arbitrary characteristic. Further developments in this direction will appear in \cite{Tyo14c, Tyo14r}.

The following geometric description of the two components of $V$, and a generalization of our example to a certain class of (not necessarily toric) complete intersection of a scroll with a
hypersurface in the projective space is due to Kristian Ranestad:

Let $K\subset \PP^n, n\geq 6$, be a threefold cone over a smooth non-degenerate projective curve $C\subset \PP^{n-2}$, and $W\subset \PP^n$ a general hypersurface of degree at least two. Then $S:=W\cap K$ is a surface with singularities on the vertex line $l$ of $K$. Let $L$ be the linear system of hyperplane sections on $S$. Then $S$ has a plane curve through every point, and the tangent planes to $S$ all intersect $l$. The Severi variety of two-nodal irreducible curves in $L$ has two components: One component corresponds to hyperplane sections of $S$  containing two tangent planes that intersect $l$ in distinct points. The other component corresponds to hyperplane sections of $S$  containing two tangent planes that intersect $l$ in the same point. The image of the nodal curves of the first component all contain $l$ in their span, while the general nodal curves of the second component do not.

To see that the pair $(X, |\CL|)$ is of the form $(S,L)$, consider the linear system of curves in $\PP^2$ given by the forms $xz^4,xyz^3,y^2z^3,xy^2z^2,x^2y^2z,xy^3z,xy^4$, and let $\widetilde{X}$ be the blowup of $\PP^2$ along the base locus of the system. Consider the map $f\:\widetilde{X}\to\PP^6$ given by the strict transforms of these forms. Then $X$ is isomorphic to the image of $f$, and $\CL$ is isomorphic to the restriction of $\CO_{\PP^6}(1)$ to $f(\widetilde{X})$. If we denote the homogeneous coordinates on $\PP^6$ by
$w_0=xz^4, w_1=xyz^3, w_2=y^2z^3, w_3=xy^2z^2, w_4=x^2y^2z, w_5=xy^3z, w_6=xy^4$
then the surface $f(\widetilde{X})$ is the complete intersection of the quadric hypersurface given by $w_2w_4=w_3^2$ with the three-dimensional cone given by
$\rank\left(
    \begin{array}{cccc}
        w_0 & w_1 & w_3 & w_5 \\
        w_1 & w_3 & w_5 & w_6 \\
    \end{array}
\right)\le 1$
i.e., the cone over the rational normal quartic curve. It is not difficult to check that $f(\widetilde{X})$ is a surface of degree eight with four isolated singularities at $[1:0:0:0:0:0:0],[0:0:1:0:0:0:0],[0:0:0:0:1:0:0],[0:0:0:0:0:0:1].$ Notice that the quadric hypersurface in this case is not general (in particular not all the singularities of the surface belong to the vertex line), nevertheless Ranestad's description of the two components of the Severi variety of two-nodal irreducible curves in $L$ is valid in this case.

{\bf Acknowledgements:} This research was partially done while the author was participating in the program {\em Tropical Geometry and Topology} at the Max-Planck Institute for Mathematics in Bonn. I am very grateful to the organizers of the program for inviting me and to MPIM Bonn for its hospitality. I would also like to thank the organizers of the {\em G\"okova Geometry and Topology conference} for inviting me to speak there.

\end{document}